\documentclass[12pt,twoside]{article}

\usepackage{amssymb}
\usepackage{amsmath}
\usepackage{bbm}
\usepackage{mathrsfs}

\makeatletter
\def\hsmash{\relax 
  \ifmmode\def\next{\mathpalette\mathhsm@sh}\else\let\next\makehsm@sh
  \fi\next}
\def\makehsm@sh#1{\setbox\z@\hbox{#1}\finhsm@sh}
\def\mathhsm@sh#1#2{\setbox\z@\hbox{$\m@th#1{#2}$}\finhsm@sh}
\def\finhsm@sh{\wd\z@\z@ \box\z@}
\makeatother

\makeatletter
\gdef\th@mychange{\normalfont\slshape
   \def\@begintheorem##1##2{\item
        [\hskip\labelsep \theorem@headerfont ##2. ##1  \,--\!--\!--\!--  ]}%
 \def\@opargbegintheorem##1##2##3{%
   \item[\hskip\labelsep \theorem@headerfont ##2. ##1\ {\upshape(}##3{\upshape)}. \,-----  ]}}
\makeatother

\RequirePackage{theorem}
\theoremstyle{mychange}

{\theorembodyfont{\rmfamily}\newtheorem{ttt}{}[section]}
{\theorembodyfont{\rmfamily}\newtheorem{remark}[ttt]{Remark.}}
{\theorembodyfont{\rmfamily}\newtheorem{rems}[ttt]{Remarks.}}
{\theorembodyfont{\rmfamily}\newtheorem{example}[ttt]{Example.}}
{\theorembodyfont{\rmfamily}\newtheorem{nota}[ttt]{Notation.}}

{\theorembodyfont{\itshape}\newtheorem{obss}[ttt]{Observations.}}
{\theorembodyfont{\itshape}}
{\theorembodyfont{\itshape}\newtheorem{theo}[ttt]{Theorem.}}
{\theorembodyfont{\itshape}}
{\theorembodyfont{\itshape}\newtheorem{lem}[ttt]{Lemma.}}

{\theorembodyfont{\rmfamily}\newtheorem{remo}[ttt]{Remark}}

{\theorembodyfont{\rmfamily}\newtheorem{prooi}[ttt]{Proof of Theorem~\ref{Mordell}.ii).}}
{\theorembodyfont{\rmfamily}\newtheorem{prooii}[ttt]{Proof of Theorem~\ref{Mordell}.iii).}}

\newcounter{abc}

\newcounter{iii}
\newenvironment{iii}{\begin{list}{\rm \roman{iii}) }{\usecounter{iii} \leftmargin=0.0pt \labelsep=0.0pt \listparindent=0.0pt \labelwidth=0.0pt \parsep=\smallskipamount \itemsep=0.0pt \topsep=0.0pt \partopsep=\smallskipamount}}{\end{list}}

\newcommand{\Pb}{\mathop{\text{\bf P}}\nolimits}
\newcommand{\Spec}{\mathop{\text{\rm Spec}}\nolimits}
\newcommand{\tr}{\mathop{\text{\rm tr}}\nolimits}
\newcommand{\Gal}{\mathop{\text{\rm Gal}}\nolimits}

\newcommand{\bbA}{{\mathbbm A}}
\newcommand{\bbF}{{\mathbbm F}}
\newcommand{\bbQ}{{\mathbbm Q}}
\newcommand{\bbR}{{\mathbbm R}}
\newcommand{\bbZ}{{\mathbbm Z}}

\newcommand{\Fp}{\bbF_{\!p}}
\newcommand{\Fq}{\bbF_{\!q}}
\newcommand{\notd}{\,\mathord{\nmid}\,}

\newcommand{\modulo}{\textup{ mod }}

\newcommand{\br}{ }
\newcommand{\brr}{, }

\def\rightend#1#2{{%
 \leavevmode\nobreak\hskip .5em plus 1fil
 \penalty600 \hskip 0pt plus -1filll
 \vadjust{}\nobreak\hskip 0pt plus 1filll%
 #1\parfillskip=#2\relax \par}}

\def\eop{\ifmmode\rule[-22pt]{0pt}{1pt}\ifinner\tag*{$\square$}\else\eqno{\square}\fi\else\rightend{$\square$}{0pt}\fi}

\renewcommand{\thefootnote}{\arabic{footnote}}
\author{J\"org Jahnel}
\date{}

\title{More cubic surfaces violating the Hasse principle}

\begin{document}
\renewcommand{\thefootnote}{\fnsymbol{footnote}}

\maketitle

\begin{abstract}
We generalize L.\,J.~Mordell's construction~\cite{Mordell} of cubic surfaces for which the Hasse principle~fails.

\centerline{\bf R\'esum\'e}
Nous g\'en\'eralisons la construction due \`a L.\,J.~Mordell~\cite{Mordell} des surfaces cubiques pour lesquelles la principe de Hasse est~fausse.
\end{abstract}

\section{Introduction and main result}

\begin{ttt}
Sir P.~Swinnerton-Dyer~\cite{SD62} was the first to construct a cubic surface
over~$\bbQ$
for which the Hasse~principle provably~fails. Swinnerton-Dyer's construction had soon been generalized by L.\,J.~Mordell~\cite{Mordell} who found two series of such~examples.\linebreak[3] The~starting points of Mordell's construction are the cubic number fields contained
in~$\bbQ(\zeta_p)$
for~$p = 7$
and~$p = 13$,~respectively.

The~failure of the Hasse~principle in Mordell's examples may be explained more conceptually by the Brauer-Manin~obstruction. This~was observed by Yu.~I.~Manin in his book~\cite[Proposition~47.4]{Ma}.
\end{ttt}

\begin{ttt}
In~this note, we will show that Mordell's construction may be~generalized to an arbitrary
prime~$p \equiv 1 \pmod 3$.
\end{ttt}

\begin{nota}
%
\begin{iii}
\item
We~denote
by~$K/\bbQ$
the unique cubic field extension contained in the cyclotomic
extension~$\bbQ (\zeta_p) / \bbQ$.
\item
We~fix the explicit generator
$\theta \in K$
given by
$\theta := \tr\nolimits_{\bbQ (\zeta_{p}) / K} (\zeta_{p} - 1)$.
More~concretely,
$\theta = -n + \sum_{i \in (\bbF_{\!p}^*)^3} \zeta_{p}^i$
for~$n := \frac{p-1}3$.
\end{iii}
\end{nota}

\begin{theo}
\label{Mordell}
Consider the cubic surface\/
$X \subset \Pb^3_\bbQ$,
given~by
$$T_3 (a_1 T_0 + d_1 T_3) (a_2 T_0 + d_2 T_3) = \prod_{i = 1}^3 \big( T_0 + \theta^{(i)} T_1 + (\theta^{(i)})^2 T_2 \big) \, .$$
Here,~$a_1, a_2, d_1, d_2$
are integers and\/
$\theta^{(i)}$
are the images
of\/~$\theta$
under\/~$\Gal (K/\bbQ)$.

\begin{iii}
\item
Then,~the reduction\/
$X_{p}$
of\/~$X$
at\/~$p$
is given~by
$T_3 (a_1 T_0 + d_1 T_3) (a_2 T_0 + d_2 T_3) = T_0^3$.
Over~the algebraic closure,
$X_{p}$~is
the union of three planes. These~are given~by
$$T_3/T_0 = s_1, \quad T_3/T_0 = s_2, \quad T_3/T_0 = s_3$$
for\/~$s_i$
the zeroes
of\/~$T (a_1 + d_1 T) (a_2 + d_2 T) - 1 = 0$.
\item
Suppose\/~$p \notd d_1 d_2$
and\/~$\gcd (d_1, d_2) = 1$.
Then,~for every\/
$(t_0 : t_1 : t_2 : t_3) \in X(\bbQ)$,
$s := (t_3/t_0 \modulo p)$
admits the property~that
$$\frac{a_1 + d_1 s}s$$
is a cube
in\/~$\bbF_{\!p}^*$.

In~particular, if\/
$\frac{a_1 + d_1 {s_i}}{s_i} \in \bbF_{\!p}^*$
is a non-cube for every\/
$i$
such that\/
$s_i \in \Fp$
then\/
$X(\bbQ) = \emptyset$.
\item
Assume~that\/~$p \notd d_1 d_2$,
that\/~$\gcd (a_1, d_1)$
and\/~$\gcd (a_2, d_2)$
contain only prime~factors which decompose
in\/~$K$,
and
that\/~$T (a_1 + d_1 T) (a_2 + d_2 T) - 1 = 0$
has at least one zero which is simple and
in\/~$\bbF_{\!p}$.
Then,~$X (\bbA_\bbQ) \neq \emptyset$.
\end{iii}
\end{theo}

\begin{ttt}
For~$p=7$
and~$13$,
we recover exactly the result of L.\,J.~Mordell.
The original example of Sir P.~Swinnerton-Dyer reappears for
$p=7$,
$a_1 = d_1 = a_2 = 1$,
and~$d_2 = 2$.
\end{ttt}

\begin{rems}
\begin{iii}
\item
$K/\bbQ$
is an abelian cubic field~extension. It~is totally ramified
at~$p$
and unramified at all other~primes.
A~prime~$q \neq p$
is decomposed
in~$K$
if and only if
$q$
is a cube
modulo~$p$.
\item
We~will
write~$\mathfrak{p}$
for the prime ideal
in~$K$
lying
above~$(p)$.
Note~that
$\mathfrak{p} = (\theta)$
by virtue of our definition
of~$\theta$.
\item
We~have
$\smash{\prod_{i = 1}^3 \!\big( T_0 + \theta^{(i)} T_1 + (\theta^{(i)})^2 T_2 \big) = N_{K/\bbQ} (T_0 + \theta T_1 + \theta^2 T_2) \, .}$
\end{iii}
\end{rems}

\section{The proofs}

\begin{obss}
\label{allg}
\begin{iii}
\item
For\/~$s$
any solution of\/
$T (a_1 + d_1 T) (a_2 + d_2 T) - 1 = 0$,
the expression\/
$\frac{a_1 + d_1 s}s$
is well defined and non-zero
in\/~$\Fp$.
\item
No\/~$\bbQ_p$-valued~point
on\/~$X$
reduces to the triple~line
``$\,T_0 = T_3 = 0$''.
\item
For~every\/~$(t_0 : t_1 : t_2 : t_3) \in X(\bbQ_p)$,
the fraction\/
$\frac{a_1 t_0 + d_1 t_3}{t_3}$
is a\/
\mbox{$p$-adic~unit.}
\end{iii}\medskip

\noindent
{\bf Proof.}
{\em
i)
By~assumption, we have
$s \neq 0$
and~$a_1 + d_1 s \neq 0$.\smallskip

\noindent
ii)
Suppose,~$(t_0 : t_1 : t_2 : t_3) \in X (\bbQ_p)$
is a point reducing to the triple~line. We~may assume
$t_0, t_1, t_2, t_3 \in \bbZ_p$
are~coprime.
Then~$\nu_p (t_0) \geq 1$
and~$\nu_p (t_3) \geq 1$
together imply that
$\smash{\nu_p \big( t_3 (a_1 t_0 + d_1 t_3) (a_2 t_0 + d_2 t_3) \big) \geq 3}$.
On~the other~hand,
$\smash{\nu_p \big( \prod_{i = 1}^3 ( t_0 + \theta^{(i)} t_1 + (\theta^{(i)})^2 t_2 ) \big) = 1}$
or~$2$
since~$t_1$
or~$t_2$
is a
unit~and~$\smash{\nu_p (\theta^{(i)}) = \frac13}$.\smallskip

\noindent
iii)
Again,~assume
$t_0, t_1, t_2, t_3 \in \bbZ_p$
to be~coprime. Assertion~ii) implies that
$t_3$~is
a
\mbox{$p$-adic~unit}.
Hence,~$\frac{a_1 t_0 + d_1 t_3}{t_3} \in \bbZ_p$.
Further,~$\smash{(\frac{a_1 t_0 + d_1 t_3}{t_3} \modulo p) = \frac{a_1 + d_1 s}s}$
for
$s := (t_3/t_0 \modulo p)$
a solution
of~$T (a_1 + d_1 T) (a_2 + d_2 T) - 1 = 0$.
}
\eop
\end{obss}

\begin{lem}
\label{good}
Let\/~$X$
be as in Theorem~\ref{Mordell}.ii) and\/
$(t_0 : t_1 : t_2 : t_3) \in X(\bbQ_\nu)$.
Further,~let\/
$\nu$
be any valuation
of\/~$\bbQ$
different
from\/~$\nu_{\!p}$
and\/~$w$
an extension
of\/~$\nu$
to\/~$K$.\smallskip

\noindent
Then,~$\frac{a_1 t_0 + d_1 t_3}{t_3} \in \bbQ_\nu^*$
is in the image of the norm
map\/~$N\colon K_w \to \bbQ_\nu$.\medskip

\noindent
{\bf Proof.}
{\em
{\em First step:}
Elementary cases.\smallskip

\noindent
If~$q$
is a prime decomposed
in~$K$
then every element
of~$\bbQ_q^*$
is a~norm. The~same applies to the infinite~prime.\medskip

\noindent
{\em Second step:}
Preparations.\smallskip

\noindent
It~remains to consider the case that
$q$
remains prime
in~$K$.
Then,~an
element~$x \in \bbQ_q^*$
is a norm if and only
if~$3 | \nu (x)$
for~$\nu := \nu_q$.

It~might happen that
$\theta$~is
not a~unit
in~$K_w$.
However,~as~$K_w / \bbQ_q$
is unramified, there exists some
$t \in \bbQ_q^*$
such that
$\underline\theta := t \theta \in K_w$
is a~unit. The~surface
$\smash{\widetilde{X}}$
given by
$$T_3 (a_1 T_0 + d_1 T_3) (a_2 T_0 + d_2 T_3) = \prod_{i = 1}^3 \big( T_0 + \underline\theta^{(i)} T_1 + (\underline\theta^{(i)})^2 T_2 \big)$$
is isomorphic
to~$X \times_{\Spec \bbQ} \Spec \bbQ_q$.
Even~more, the map
$$\textstyle \iota\colon X \times_{\Spec \bbQ} \Spec \bbQ_q \to \widetilde{X}, \quad (t_0 : t_1 : t_2 : t_3) \mapsto (t_0 : \frac{t_1}t : \frac{t_2}{t^2} : t_3)$$
is an isomorphism which leaves the rational
function~$\smash{\frac{a_1 T_0 + d_1 T_3}{T_3}}$~unchanged.
Hence,~we may assume without restriction
that~$\theta \in K_w$
is a~unit.\medskip

\noindent
{\em Third step:}
The~case that
$\theta$
is a~unit.\smallskip

\noindent
Assume~that
$t_0, t_1, t_2, t_3 \in \bbZ_q$
are~coprime.
If~$\nu \big( t_3 (a_1 t_0 + d_1 t_3) (a_2 t_0 + d_2 t_3) \big) = 0$
then
$\smash{\frac{a_1 t_0 + d_1 t_3}{t_3}}$
is clearly a~norm. Otherwise,~we
have~$\smash{\nu \big( \prod_{i = 1}^3 ( t_0 + \theta^{(i)} t_1 + (\theta^{(i)})^2 t_2 ) \big) > 0}$.
This~implies that
$\nu (t_0)$,~$\nu (t_1)$,~$\nu (t_2) > 0$.
Consequently,
$t_3$~must
be a~unit.

From~the equation
of~$X$,
we
deduce~$\nu (d_1 d_2) > 0$.
If~$\nu (d_2) > 0$
then, according to the assumption,
$d_1$~is
a~unit. This~shows
$\nu (a_1 t_0 + d_1 t_3) = 0$
from which the assertion~follows.

Thus,~assume~$\nu (d_1) > 0$.
Then,~$d_2$
is a~unit and, therefore,
$\nu (a_2 t_0 + d_2 t_3) = 0$.
Further,~we note~that
$$3 \,|\, \nu \Big(\! \prod_{i = 1}^3 ( t_0 + \theta^{(i)} t_1 + (\theta^{(i)})^2 t_2 ) \!\Big)$$
since the product is a~norm. By~consequence,
$3 | \nu \big( t_3 (a_1 t_0 + d_1 t_3) (a_2 t_0 + d_2 t_3) \big)$.
Altogether,~we see
that~$3 | \nu (a_1 t_0 + d_1 t_3)$
and~$\smash{3 | \nu \big( \frac{a_1 t_0 + d_1 t_3}{t_3} \big)}$.
The~assertion~follows.
}
\eop
\end{lem}

\begin{prooi}
According~to Lemma~\ref{good},
$\frac{a_1 t_0 + d_1 t_3}{t_3} \in \bbQ^*$
is a local norm at every prime
except~$p$.
Global~class field theory~\cite[Theorem~5.1 together with 6.3]{Ta} shows that it must necessarily be a norm at that prime,~too.

$(p) = {\mathfrak p}^3$
is a totally ramified~prime.
A~$p$-adic
unit~$u$
is a norm if and only if
$\overline{u} := (u \modulo p)$
is a cube
in~$\bbF_{\!p}^*$.
By~Observation~\ref{allg}.iii),
$\smash{\frac{a_1 t_0 + d_1 t_3}{t_3}}$
is automatically
a~$p$-adic~unit.
As~$\smash{\frac{a_1 + d_1 s}s = (\frac{a_1 t_0 + d_1 t_3}{t_3} \modulo p)}$,
this is exactly the~assertion.
\eop
\end{prooi}

\begin{prooii}
We~have to show that
$X(\bbQ_\nu) \neq \emptyset$
for every valuation
of~$\bbQ$.
$X(\bbR) \neq \emptyset$
is~obvious.
For~a
prime
number~$q$,
in~order to prove
$X(\bbQ_q) \neq \emptyset$,
we use Hensel's~lemma.
It~is sufficient to verify that the reduction
$X_q$
has a smooth
$\Fq$-valued~point.
Thereby,~we may replace
$X$
by a
\mbox{$\bbQ_q$-scheme}~$\smash{\widetilde{X}}$
isomorphic
to~$X \times_{\Spec \bbQ} \Spec \bbQ_q$.\medskip

\noindent
{\em Case~1:}
$q = p$.\smallskip

\noindent
Then,~the
reduction~$X_{p}$
is the union of three planes meeting in the~line given
by~$T_0 = T_3 = 0$.
By~assumption, one of the planes appears with multiplicity one and is defined
over~$\Fp$.
It~contains
$p^2$
smooth~points.\medskip

\noindent
{\em Case~2:}
$q \neq p$.\smallskip

\noindent
Assume~without restriction that
$\theta$~is a
\mbox{$w$-adic}~unit.
There~are two~subcases.\smallskip

\noindent
a)
$q \notd d_1 d_2$.\hspace{2.5mm}
It~suffices to show that there is a smooth
$\Fq$-valued~point
on the
intersection~$X_q^\prime$
of~$X_q$
with the hyperplane~``$T_0 = 0$''.
This~curve is given~by
$$\overline{d}_1 \overline{d}_2 T_3^3 = \overline\theta^{(1)} \overline\theta^{(2)} \overline\theta^{(3)} \prod_{i = 1}^3 (T_1 + \overline\theta^{(i)} T_2) \, .$$
If~$q \neq 3$
then this equation defines a smooth genus one~curve. It~has an
$\Fq$-valued
point by Hasse's~bound.

If~$q = 3$
then the projection
$X_q^\prime \to \Pb^1$
given by
$(T_1 : T_2 : T_3) \mapsto (T_1 : T_2)$
is one-to-one
on~$\Fq$-valued~points.
At~least one of them is smooth since
$\smash{\prod_{i=1}^3 (T + \overline\theta^{(i)})}$
is a separable~polynomial.\smallskip

\noindent
b)
$q | d_1 d_2$.\hspace{2mm}
$X_q^\prime := X_q \;\cap\;$``$T_0 = 0$''
is given by
$\smash{0 = \overline\theta^{(1)} \overline\theta^{(2)} \overline\theta^{(3)} \prod_{i = 1}^3 (T_1 + \overline\theta^{(i)} T_2)}$.
In~particular,
$x = (0 : 0 : 0 : 1) \in X_q (\Fq)$.
We~may assume that
$x$
is~singular.

Then,~$X_q$
is given as
$Q (T_0, T_1, T_2) T_3 + K (T_0, T_1, T_2) = 0$
for~$Q$
a quadratic form
and~$K$
a cubic form.
If~$Q \not\equiv 0$
then there is
an~$\Fq$-rational
line~$\ell$
through~$x$
such that
$Q |_\ell \neq 0$.
Hence,~$\ell$
meets~$X_q$
twice
in~$x$
and once in another
$\Fq$-valued~point
which is~smooth.

Otherwise,~$(F \modulo q)$
does not depend
on~$T_3$.
I.e.,~the left hand side of the equation
of~$X$
vanishes
modulo~$q$.
This~means that one of the factors must~vanish. We~have, say,
$a_1 \equiv d_1 \equiv 0 \pmod q$.
Then,~by assumption,
$q$~decomposes
completely
in~$K$.
At~such a prime,
$X_q^\prime$~is
the union of three lines which are defined
over~$\Fq$,
different from each~other, and meet in one~point. There~are plenty of smooth points
on~$X_q^\prime$.
\eop
\end{prooii}

\section{Examples}

\begin{example}
For~$p = 19$,
a counterexample to the Hasse~principle is given~by
\begin{eqnarray*}
T_3 (19 T_0 + 5 T_3) (19 T_0 + 4 T_3) & = & \prod_{i = 1}^3 \big( T_0 + \theta^{(i)} T_1 + (\theta^{(i)})^2 T_2 \big) \\
 & = & T_0^3 - 19 T_0^2 T_1 + 133 T_0^2 T_2 + 114 T_0 T_1^2 \\
 &   & \hspace{1.6cm}{} - 1\,539 T_0 T_1 T_2 + 5\,054 T_0 T_2^2 - 209 T_1^3 \\
 &   & \hspace{2.4cm}{} + 3\,971 T_1^2 T_2 - 23\,826 T_1 T_2^2 + 43\,681 T_2^3 \, .
\end{eqnarray*}
Indeed,~in~$\bbF_{\!19}$,
the cubic~equation
$T^3 - 1 = 0$
has the three solutions
$1$,
$7$,
and~$11$.
However,~in any case
$\frac{a_1 + d_1 s}s = 5$
which is a non-cube.
\end{example}

\begin{example}
\label{Angabe}
Put~$p = 19$.
Consider~the cubic
surface~$X$
given~by
$$T_3 (T_0 + T_3) (12 T_0 + T_3) = \prod_{i = 1}^3 \big( T_0 + \theta^{(i)} T_1 + (\theta^{(i)})^2 T_2 \big) \, .$$
Then,~$X (\bbA_\bbQ) \neq \emptyset$
but~$X (\bbQ) = \emptyset$.
$X$~violates
the Hasse~principle.

Indeed,~in~$\bbF_{\!19}$,
the cubic~equation
$T (1 + T) (12 + T) - 1 = 0$
has the three solutions
$12$,
$15$,
and~$17$.
However,~in~$\bbF_{\!19}$,
$13/12 = 9$,
$16/15 = 15$,
and~$18/17 = 10$
which are three non-cubes.
\end{example}

\begin{example}
For~$p = 19$,
consider~the cubic
surface~$X$
given~by
$$T_3 (T_0 + T_3) (2 T_0 + T_3) = \prod_{i = 1}^3 \big( T_0 + \theta^{(i)} T_1 + (\theta^{(i)})^2 T_2 \big) \, .$$
Then,~for~$X$,
the Hasse principle~fails.

Indeed,~in~$\bbF_{\!19}$,
the cubic~equation
$T (1 + T) (2 + T) - 1 = 0$
has~$T = 5$
as its only~solution. The~two other solutions are conjugate to each other
in~$\bbF_{\!19^2}$.
However,~in~$\bbF_{\!19}$,
$6/5 = 5$
is a non-cube.
\end{example}

\begin{example}
\label{obsweak}
Put~$p = 19$
and consider~the cubic
surface~$X$
given~by
$$T_3 (T_0 + T_3) (6 T_0 + T_3) = \prod_{i = 1}^3 \big( T_0 + \theta^{(i)} T_1 + (\theta^{(i)})^2 T_2 \big) \, .$$
There~are
\mbox{$\bbQ$-rational}
points
on~$X$
but weak approximation~fails.

Indeed,~in~$\bbF_{\!19}$,
the cubic~equation
$T (1 + T) (6 + T) - 1 = 0$
has the three solutions
$8$,
$9$,
and~$14$.
However,~in~$\bbF_{\!19}$,
$10/9 = 18$
is a cube while
$9/8 = 13$
and~$15/14 = 16$
are non-cubes.
The~smallest
$\bbQ$-rational
point
on~$X$
is~$(14 : 15 : 2 : (-7))$.
Observe~that, in fact,
$T_3 / T_0 = -7/14 \equiv 9 \pmod {19}$.
\end{example}

\begin{remark}
From~each of the examples given, by adding multiplies
of~$p$
to the coefficients
$a_1$,
$d_1$,
$a_2$,
and~$d_2$,
a family of surfaces arises which are of similar~nature.
\end{remark}

\begin{remo}[{\rm Lattice basis reduction}{}]
The~norm form in the
$p = 19$~examples
produces coefficients which are rather~large. An~equivalent form with smaller coefficients may be obtained using lattice basis reduction. In~its simplest form, this means the~following.

For~the
rank~$2$
lattice
in~$\bbR^3$,
generated by the
vectors~$v_1 := (\theta^{(1)}, \theta^{(2)}, \theta^{(3)})$
and
$v_2 := \big( (\theta^{(1)})^2, (\theta^{(2)})^2, (\theta^{(3)})^2 \big)$,
in~fact
$\{ v_1, v_2 + 7 v_1 \}$
is a reduced~basis. Therefore,~the substitution
$T_1^\prime := T_1 - 7 T_2$
simplifies the norm~form. Actually,~we~find
\begin{eqnarray*}
\smash{\prod_{i = 1}^3 \!\big( T_0 + \theta^{(i)} T_1 + (\theta^{(i)})^2 T_2 \big)} & = & \!T_0^3 \!-\! 19 T_0^2 T_1^\prime \!+\! 114 T_0 T_1^\prime {}^2 \!+\! 57 T_0 T_1^\prime T_2 \!-\! 133 T_0 T_2^2 \\
 &  & \hspace{1.5cm}{} - 209 T_1^\prime {}^3 - 418 T_1^\prime {}^2 T_2 + 1045 T_1^\prime T_2^2 - 209 T_2^3 \, .
\end{eqnarray*}
\end{remo}

\end{document}